\input amstex.tex
\input amsppt.sty
\hsize=173truemm
\vsize=238truemm
\voffset=-.3truecm
\hoffset=-1.2truecm
\def\({\left(}
\def\){\right)}

\topmatter
\title
%Спектральные оценки оператора Лапласа дискретной группы Гейзенберга
Spectral estimations for Laplace operator for the discrete Heisenberg
group
\endtitle
\author K.~Kokhas, A.~Suvorov \endauthor

\abstract\nofrills{}
We estimate the spectral measure of Laplace operator
$\Delta=\frac{1}{4}(x+x^{-1}+y+y^{-1})$ for the discrete Heisenberg
group with generators $x$ and $y$ in vicinity of the unity.
\endabstract
\endtopmatter
\bigskip

Let $H$ be the discrete Heisenberg group
$$
H=\Bigl\{
\pmatrix
1&k&m\cr
0&1&\ell\cr
0&0&1
\endpmatrix,\quad
k, \ell, m \in {\Bbb Z}\Bigr\}\,.
$$
We will call analogous matrix group with elements from
${\Bbb Z}/n{\Bbb Z}$ a finite Heisenberg group and denote it $H_{n}$.
The generators of $H$ and $H_n$ will be denoted
by the same letters:
$$
x=\pmatrix
1&1&0\cr
0&1&0\cr
0&0&1
\endpmatrix,\quad
y=\pmatrix
1&0&0\cr
0&1&1\cr
0&0&1
\endpmatrix,\quad
z=\pmatrix
1&0&1\cr
0&1&0\cr
0&0&1
\endpmatrix\,. \eqno(1)
$$
The element $\Delta$ of the group algebra for~$H$,
$\Delta=\frac{1}{4}(x+x^{-1}+y+y^{-1})$,
is the Laplace operator corresponding to the system of generators $(x, y)$.
Analogously,
$\Delta_n=\frac{1}{4}(x+x^{-1}+y+y^{-1})$ is the Laplace operator for the
finite Heisenberg group~$H_n$.
These operators can also be defined as transition operators for random walk
on the groups.
We will deal with spectra of these operators in the regular
representations of the corresponding groups.

It is easy to demonstrate [BVZ] that the spectrum of
$\Delta$ is the interval $[-1,1]$.
Let $E_A$, $A\subset [-1,1]$ be a family of spectral projectors for
$\Delta$ and
$\mu{A}=(E_A\delta_e,\delta_e)$ be the corresponding spectral measure.
Here $\delta_e\in L_2(H)$ is the characteristic function of the unit element
of the group~$H$.
We will estimate the value
$\mu\bigl([-1,-1+t]\cup [1-t,1]\bigr)$  for $t\to 0$.
More precisely we will prove the inequality
$\mu\bigl([-1,-1+t]\cup [1-t,1]\bigr)\geq \text{const}\,t^{2+\alpha}$.

The paper is organized as follows.
The first section contains description of the representations of the finite
Heisenberg group.
In the second section we demonstrate that all eigenvalues of the operator
$\Delta_n$ in multidimensional representations are less than
$1-O(\frac{1}{n})$.
In the third section we give a combinatorial realization for the
characteristic polynomials of $\Delta_n$ in the
multidimensional representations.
And in the fourth section we prove the estimates for the spectral measure of
$\Delta$, using those for the finite group.

The work was partially supported by the grant RFFI 96-15-96060.

\subhead
1. Representations of the finite Heisenberg group $H_n$
\endsubhead
%\bigskip

We will consider only the simplest case: $n$ is a prime number.
In this case the group  $H_n$ has $n^2$ one-dimensional representations
$T_{\alpha,\beta}$, $\alpha$, $\beta=1$, 2, \dots, $n$: \
$$
T_{\alpha,\beta}(x)=e^{\frac{2\pi i}{n} \alpha}, \quad
T_{\alpha,\beta}(y)=e^{\frac{2\pi i}{n} \beta},  \quad
T_{\alpha,\beta}(z)=1\,.
$$
Besides there are $n-1$ irreducible representations $T_q$ of dimension~$n$.
These representations can be obtained by induction from the representations
$\rho_q$ of the abelian subgroup generated by the elements $y$ and $z$,
where  $\rho_{q}(y)=1$, $\rho_q(z)=e^{\frac{2\pi i}{n} q}$, $q=1$, \dots,
$n-1$.
The representation $T_q$ can be described in the space ${\Bbb C}^n$
with the basis $u_1$, $u_2$, \dots, $u_n$ by the formulae
$$
T_q(x)\,u_j = u_{j+1},\quad
T_q(y)\,u_j = e^{\frac{2\pi i}{n} qj}u_j,\quad
T_q(z)\,u_j = e^{\frac{2\pi i}{n} q}u_j\,.
$$

\proclaim{Proposition 1.1}  $T_{\alpha,\beta}$ and $T_q$ is a complete
set of irreducible unitary representations of the group $H_n$.
\endproclaim

\smallskip
Indeed, irreducibility and nonequivalence of the representations $T_q$
can be verified by computing scalar products of the corresponding
characters. Computing squares of the dimensions of the representations
$T_{\alpha,\beta}$ and $T_q$ we can see that the set is complete.
\smallskip

Let  $\Delta_n=\frac{1}{4}(x+y+x^{-1}+y^{-1})\in {\Bbb C}H_n$ be
the Laplace operator corresponding to the generators~(1).
Denote $\tilde{\Delta}_n=4\Delta_n$. In order to
find the spectrum of $\tilde{\Delta}_n$ in the regular representation
it suffices to find its spectra in all irreducible representations.
In the case of the one-dimensional representations $T_{\alpha,\beta}$
the operator is
$$
T_{\alpha,\beta}(\tilde{\Delta}_n)=
 2\cos\frac{2\pi\alpha }{n}+ 2\cos\frac{2\pi\beta}{n}\,.
$$
For multidimensional representations $T_q$ the operator has a more complicated
form
$$
T_q(\tilde{\Delta}_n)=
\pmatrix
2\cos\frac{0\cdot2q\pi}{n}&1&0&0&\dots&1\cr
1&2\cos\frac{1\cdot2q\pi}{n}&1&0&\dots&0\cr
0&1&2\cos\frac{2\cdot2q\pi}{n}&1&\dots&0\cr
%\hdotsfor6\cr
\vspace{3pt}
\vdots&\vdots&\vdots&\vdots&\ddots&\vdots\cr
\vspace{6pt}
0&0&\dots&1&2\cos\frac{(n-2)\cdot2q\pi}{n}&1\cr
1&0&\dots&0&1&2\cos\frac{(n-1)\cdot2q\pi}{n}\cr
\endpmatrix
$$
Unfortunately the spectra of these matrices cannot be computed directly.

Such matrices appear in investigations of the Harper operator
and the Mathieu equation (see e.g. \hbox{[AMS, BVZ]}).
Consider the Harper operator in the space $L^2({\Bbb Z})$
$$
u_n \mapsto u_{n+1} + u_{n-1} + 2\cos(2\pi n\alpha) u_n\,.
$$
Specialists in the high energy physics know very well the so called
Hofstadter butterfly [H], depicted in the figure~1 (we have taken this
picture from~[IMOS]). Here the horizontal axis corresponds to the parameter
$\alpha$, and the vertical axis~--- to the spectrum.
If $\alpha =\frac{q}n$ is a rational number (irreducible fraction),
then the spectrum consists of  $n$ intervals of the form
$f^{-1}([-4, 0])$, where $f$ is the characteristic polynomial of the matrix
$T_q(\Delta_n)$ (if $n$ is even the $\frac{n}2$-th and $(\frac{n}{2}+1)$-th
intervals have the common vertex~0).
If  $\alpha$ is irrational, the spectrum is a Cantor set.

The endpoints of the intervals form curves well seen in the picture.
The curve marked by the arrow has a slope $c\ne 0$ and it implies that
$\rho_n\sim 4-\frac{c}{n}$, where $\rho_n$ is the spectral radius
of the Harper operator for $\alpha = \frac{1}{n}$.
%Поскольку характеристический
%многочлен сильно осциллирует на промежутке $[-1,1]$,
In fact $\rho_n =\lambda_n$, where $\lambda_n$ is the maximal eigenvalue of
the matrix $T_1(\Delta_n)$. In the next section we will show that
$4-\lambda_n= O(\frac 1n)$, which agrees with the form of the marked curve.

Numerical computations allow us to conjecture that
$\lambda_n = 4 - \frac{2\pi}{n} + o(\frac{1}{n})$ when $n\to +\infty$.
For instance, for $n=1000$ the value of  $(4-\lambda_n) n$ is approximately
6.28 (the next digit does not fit $2\pi$; this has been computed by means of
the Mathematica 3.0).

\bigbreak
\subhead 2. Estimation of the maximal eigenvalue of the operator
$\tilde{\Delta}_n$ in the representation $T_q$.
\endsubhead

\proclaim{Notations} \rm Matrices are denoted by ``calligraphic''
letters~--- $\Cal A$, $\Cal B$ etc.  The maximal eigenvalue of ${\Cal A}$
is denoted~$\lambda_{\Cal A}$. Let
${\Cal P}=(p_{ij})$ and ${\Cal Q}=(q_{ij})$ be  real $n\times n$ matrices.
We write ${\Cal P}\leq {\Cal Q}$, if $p_{ij}\leq q_{ij}$ for all  $i$,~$j$.
\endproclaim

Let ${\Cal A}_n$ be $n\times n$ matrix defining our operator
$\tilde{\Delta}_n$ in the representation~$T_1$:
$$
{\Cal A}_n=
\pmatrix
2\cos\frac{0\cdot2\pi}{n}&1&0&0&\dots&1\cr
1&2\cos\frac{1\cdot2\pi}{n}&1&0&\dots&0\cr
0&1&2\cos\frac{2\cdot2\pi}{n}&1&\dots&0\cr
%\hdotsfor6\cr
\vspace{3pt}
\vdots&\vdots&\vdots&\vdots&\ddots&\vdots\cr
\vspace{6pt}
0&0&\dots&1&2\cos\frac{(n-2)\cdot2\pi}{n}&1\cr
1&0&\dots&0&1&2\cos\frac{(n-1)\cdot2\pi}{n}\cr
\endpmatrix
$$
This is a symmetric matrix, it has a real spectrum.
Let $\lambda _n$ be its maximal eigenvalue.
It is convenient to deal with the matrix
$\tilde{\Cal A}_n=2{\Cal E}_n+{\Cal A}_n$
(${\Cal E}_n$ is the $n\times n$ identity matrix),
with non-negative entries, and the spectrum is shifted by 2 with respect to
the spectrum of ${\Cal A}_n$, and the set of eigenvectors coincides with that
of~${\Cal A}_n$. In particular, $\tilde{\lambda}_n=\lambda_n+2$ where
${\tilde{\lambda}_n}$  is the maximal eigenvalue of~$\tilde{\Cal A}_n$.
According to the Perron-Frobenius theorem the maximal eigenvalue of
an irredundant non-negative matrix has the multiplicity 1, hence
$\lambda_n$ is an eigenvalue of multiplicity 1.

We need also the following consequence of the Perron-Frobenius theorem:

\proclaim{Proposition 2.1} Let ${\Cal P}$ and ${\Cal Q}$ be $n\times n$
matrices with non-negative entries,
${\Cal P}\leq {\Cal Q}$, then $\lambda_{\Cal P}\leq\lambda_{\Cal Q}$.
\endproclaim

Let us demonstrate that $4-\lambda_n= O(\frac 1n)$ when $n\to+\infty$.
We choose the numbers $\sqrt{n/2}$, $37/15$ etc in the proofs of the
next two lemmas in such a way that we obtain good approximations of exact
constants.

\proclaim{Lemma 2.2} $4-\frac{40}{n}<\lambda_n<4-\frac{2}{n}$\,.
\endproclaim

\demo{Proof} 1) The left inequality.

Since it does not affect the final result, we will treat the expression
$\sqrt {n/2}$ as even integer. Let $\botsmash{{\Cal C}_{\sqrt {n/2}}}$ be
tridiagonal $\sqrt {n/2}\, \times\sqrt {n/2}$ matrix of the form
$$
{\Cal C}_{\sqrt {n/2}}=
\pmatrix
0&1&0&\dots&0&0\\
1&0&1&\dots&0&0\\
0&1&0&\dots&0&0\\
\vdots&\vdots&\vdots&\ddots&\vdots&\vdots\\
0&0&0&\dots&1&0\\
\endpmatrix\,.
$$
It is well known (see e.g.~[GHJ]) that the spectrum of
${\Cal C}_{\sqrt {n/2}}$ is the set
$\bigl\{2\cos\frac{j\pi}{\strut\sqrt {n/2} +1},
\ 1\leq j\leq \sqrt {n/2}\,\bigr\}$.
Let ${\Cal B}_n$ be $n\times n$ matrix having the following block form:
$$
{\Cal B}_{n}=
\pmatrix
\bigl(4-\frac{20}{n}\bigr){\Cal E}_{\sqrt {n/2}} + {\Cal C}_{\sqrt {n/2}}&
        \bold 0&\bold 0\\
\bold 0&\bold 0&\bold 0\\
\bold 0&\bold 0&
      \bigl(4-\frac{20}{n}\bigr){\Cal E}_{\sqrt {n/2}} + {\Cal C}_{\sqrt {n/2}} \\
\endpmatrix\,.
$$
Then $\tilde{\Cal A}_n \geq {\Cal B}_n$ (since
$2\cos\frac{2\pi k}{n} \geq 2-\frac{2\pi^2}{n} >2-\frac{20}{n}$
for $|k|\leq\sqrt{n/2}$).
By the proposition  2.1 we get
$$
\lambda_n\geq \lambda_{{\Cal B}_n}-2=
2-\frac{20}{n}+2\cos\frac{\pi}{\sqrt {n/2} +1} >4-\frac{40}{n}\,.
$$

\smallskip
2) The right inequality. To simplify notations we will treat
the expression $\sqrt n /4$ as integer.

Let  ${\Cal D}_n$ be  ``generalized''-tridiagonal matrix, containing
on the main diagonal the numbers $\frac{37}{15n}$ ($\sqrt n/4$ times), then
$n+1-2\sqrt n/4$ zeros, and then again the numbers $\frac{37}{15n}$
($\sqrt n/4-1$ times):
$$
{\Cal D}_n=
\pmatrix
\frac{37}{15n}&1& & & & & & & &1 \\
1&\frac{37}{15n}&1 \\
  &\hdotsfor3&  \\
  & &1&\frac{37}{15n}&1&       \\
  & & &1&0&1 \\
  & & & &1&0&1 \\
  & & & & &\hdotsfor3& \\
  & & & & &1&0&1   \\
  & & & & & &1&\frac{37}{15n}&1  \\
  & & & & & & &\hdotsfor3     \\
1 & & & & & & & & 1&\frac{37}{15n}  \\
\endpmatrix
$$
Then the inequality
$$
\tilde{\Cal A}_n\leq \bigl(4-\frac{37}{15n}\bigr){\Cal E}_n+{\Cal D_n}\,,
\tag2
$$
holds, because for large $n$ and $\sqrt n/4\leq k \leq n-\sqrt n/4$ the
estimation $2\cos\frac{2\pi k}{n}\leq 2 -\frac{\pi^2}{4n}+O(\frac{1}{n^4}) <
2-\frac{37}{15n}$ is valid since $\pi^2/4\approx 2.467>37/15=2.466$.

Let us find an upper bound for the maximal eigenvalue $\lambda_{{\Cal D}_n}$.
It is clear, that, $2<\lambda_{{\Cal D}_n}< 2+\frac{37}{15n}$.
Let $\lambda_{{\Cal D}_n}=2+\alpha$, where $0<\alpha <\frac{37}{15n}$.
We will obtain below more precise estimation
$\alpha < \frac{5}{11n}=\frac{0.454}n$.
Then the statement of the lemma will be a direct consequence of the
inequality~(2).

We will give a proof from the contrary.
Assume that $\alpha >\frac{5}{11n}$ (the only place where we use this
assumption is the 4th step of the proof).
Using this assumption in technical calculations,  we will obtain the
inequality $\alpha <\frac{5}{11n}$. This contradiction implies that actually
$\alpha <\frac{5}{11n}$.

Let $\bold v =(v_1, v_2, \dots v_n)$ be the eigenvector corresponding to
the eigenvalue $\lambda_{{\Cal D}_n}$:
$$
{\Cal D}_n \bold v=(2+\alpha)\bold v\,.\tag 3
$$
We define  $v_k$ for all  $k\in{\Bbb Z}$ by the periodicity.
We normalize the vector $\bold v$ by the condition $\max\limits_{k} v_k=1$.
Recall that all the $v_i$ are non-negative.

Let us list some properties of  $\bold v$.

\item {1.} $\bold v$ is an eigenvector of multiplicity~1. Therefore
it is symmetric: $v_k=v_{n+2-k}$ (here $v_{n+1}=v_1$).

\item {2.} $\max\limits_{k}v_k=v_1=1$, $\min\limits_k v_k=v_{n/2}$.\hfill\break
Indeed, for $-\sqrt n/4+1\leq k\leq \sqrt n/4$ it follows from relation (3)
that
$$
v_{k-1}+v_{k+1}=\bigl(2-\frac{37}{15n}+\alpha\bigr)v_k\,,\tag4
$$
in other words $v_{k-1}+v_{k+1}<2v_k$, and $v_k$ cannot be a local minimum.
Analogously, for $\sqrt n/4<k\leq n-\sqrt n/4$
$$
v_{k-1}+v_{k+1}=(2+\alpha)v_k\,,\tag5
$$
hence $v_{k-1}+v_{k+1}>2v_k$, and $v_k$ cannot be a local maximum.
Now the statement follows from the symmetry of~$\bold v$.

\item{3.} $v_{\sqrt n/4}=v_{n-\sqrt n/4+2} > \frac{12}{13}$
for sufficiently large $n$. \hfill\break
Indeed, from the equality (4) we get
$$
v_k-v_{k+1}=\Bigl(\frac{37}{15n}-\alpha\Bigr)v_k+(v_{k-1}-v_k)\leq
\frac{37}{15n}+(v_{k-1}-v_k)\,.
$$
Summing up these relations, we obtain
$$
v_k-v_{k+1}\leq \frac{37}{15n}(k-1) + (v_1-v_2)\,,
$$
therefore
$$
v_1-v_{\sqrt n/4}\leq
\frac{37}{15n}\Bigl(1+2+\dots +(\sqrt n/4-1)\Bigr)+
                      \bigl(\sqrt n/4-1\bigr)(v_1-v_2)\,.
$$
Finally, notice that the first term does not exceed
$\frac{1}{32}\cdot\frac{37}{15}=\frac{37}{480}< \frac{1}{13}$,
and the second term is less than $\frac{37}{120\sqrt n}$
(it follows from (4) for  $k=1$ and from the symmetry of~$\bold v$).
Thus, $v_{\sqrt n/4}>\frac{12}{13}$ for large~$n$.

\item{4.} For large  $n$ \ \
$v_{\sqrt n/4+1}+v_{\sqrt n/4+2}+\dots +v_{n-\sqrt n/4}>
\frac{24}{13\sqrt\alpha}$. \hfill\break
Indeed, solving homogeneous difference equation~(5) and taking into account
the symmetry of~$\bold v$, we deduce the formulae
$$
v_{n/2+k}=C_0\(\Bigl(1+\alpha/2+\sqrt{\alpha +\alpha^2/4}\,\Bigr)^k+
\Bigl(1+\alpha/2-\sqrt{\alpha +\alpha^2/4}\,\Bigr)^k\)\quad
(-n/2+\sqrt n/4\leq k\leq n/2-\sqrt n/4)\,.
$$
So, we need only to estimate the sum of geometric progression.
Denote for brewity  $N=n/2-\sqrt n/4$,
$x=1+\alpha/2+\sqrt{\alpha +\alpha^2/4}$,
then  $x^{-1}=1+\alpha/2-\sqrt{\alpha +\alpha^2/4}$,
$$
v_{\sqrt n/4+1}+v_{\sqrt n/4+2}+\dots +v_{n-\sqrt n/4} =
C_0\sum\limits_{k=-N}^N (x^k+x^{-k}) =
2C_0\frac{x^{N+1/2}-x^{-(N+1/2)}}{x^{1/2}-x^{-1/2}}\,.\tag6
$$
By the assumption, $\frac{5}{11n}<\alpha <\frac{37}{15n}$,
$x^N=e^{\hbox{const}\sqrt n}$,
$x^{-N}=e^{-\hbox{const}\sqrt n}$, therefore
$x^{N+1/2}-x^{-(N+1/2)}\geq(1-\varepsilon)x^{N+1/2}$ for any positive
$\varepsilon $, if $n$ is sufficiently large.
Notice also that
$$
\frac{1}{x^{1/2}-x^{-1/2}}=\frac{x^{1/2}+x^{-1/2}}{x-x^{-1}}\geq
\frac{1}{\sqrt{\alpha +\alpha^2/4}}\geq\frac{1-\varepsilon}{\sqrt\alpha}\,.
$$
Now we can obtain a lower bound for the r.h.s.\ in the equality~(6):
$$
\sum\limits_{k=-N}^N v_{n/2+k}=
2C_0\frac{x^{N+1/2}-x^{-(N+1/2)}}{x^{1/2}-x^{-1/2}}\geq
(1-\varepsilon)^2\frac{2C_0x^{N}x^{1/2}}{\sqrt\alpha}=
(1-\varepsilon)^2\frac{2v_{\sqrt n/4}}{\sqrt\alpha} >
\frac{24}{13\sqrt\alpha}\,,
$$
($\varepsilon$ should be chosen sufficiently small in order to the last
inequality was correct according to the 3rd step of the proof).

\smallskip
To complete the proof of lemma 2, sum up equalities  (4) and~(5)
for all  $k$:
$$
\frac{37}{15n}\sum_{k=-\sqrt n/4+1}^{\sqrt n/4}v_k=
\alpha\sum_{k=1}^{n}v_k
$$
Due to normalization $0\leq v_i\leq 1$ and the inequality given in the 4th
step, we easily estimate both sides
$$
2\cdot\frac{37}{15n}\sqrt n/4>
\frac{37}{15n}\sum_{k=-\sqrt n/4+1}^{\sqrt n/4}v_k\geq
\alpha\!\!\sum_{k=\sqrt n/4+1}^{n-\sqrt n/4}v_k >
\frac{24}{13}\sqrt\alpha\,.
$$
Thus,
$$
\alpha\leq
\Bigl(\frac{37\cdot2\cdot13}{24\cdot15\cdot4}\Bigr)^2\cdot\frac{1}{n}<
\frac{4}{11n}\,.
$$
\hbox{}\hfill $\square$ Lemma 2.2 is proven.
\enddemo

For the maximal eigenvalues of the operator $\tilde{\Delta}_n$
in the other representations $T_k$ we will give only an upper bound.
It is worse than that of lemma  2.2, however we believe that
$\lambda_n$ in lemma 2.2 is the maximal eigenvalue among the eigenvalues
of matrices $T_q(\Delta_n)$ for all~$q$.
Let $\mu_n$ be the maximal eigenvalue of the operator $\tilde{\Delta}_n$
in any of representations~$T_q$.

\proclaim{Lemma 2.3} $\mu_n \leq 4-\frac{3}{5n}$.
\endproclaim

\demo{Proof} Similarly to the proof of the lemma 2.2, we choose a matrix
which majorizes
%????????????????
$T_q(\tilde{\Delta}_n)$. For this purpose we replace $\sqrt n/3$
largest expressions of the form $2\cos\frac{2\pi qj}{n}$ at the diagonal,
with 2 and all the other diagonal entries with $2-1/n$.
Then we have the following inequality
$T_q(\tilde{\Delta}_n)+2{\Cal E}_n\leq
(4-\frac{1}{n}){\Cal E}_n + \tilde{\Cal D}_n$,
since $\cos\frac{2\pi\sqrt{n}/6}{n}\leq 2-\frac{\pi^2}{9n}+O(\frac{1}{n^4})
<2-\frac{1}{n}$ for large~$n$.
Here $\tilde{\Cal D}_n$ is a matrix analogous to  ${\Cal D}_n$,
having $\sqrt n/3$ numbers $1/n$ at the diagonal in the places corresponding
to large cosines, and zeros in the other places.
Let $\bold v =(v_1, v_2, \dots v_n)$ be the eigenvector for the number
$\lambda_{\tilde{\Cal D}_n}$, normalized by condition $\max\limits_{k}v_k=1$.
As above, denote $\lambda_{\tilde{\Cal D}_n}=2+\alpha$.
To prove the lemma it suffices to verify that $\alpha\leq \frac{2}{5n}$.

Components of vector $\bold v$ obey the relations
$$
v_{k-1}+v_{k+1}=\bigl(2-\frac{1}{n}+\alpha\bigr)v_k
\hbox{\ \ \ or \ \ \ }
v_{k-1}+v_{k+1}=(2+\alpha)v_k\,.
$$
Multiply each relation by $v_k$ and sum over all $k$:
$$
2\sum\limits_{k=1}^n v_k v_{k+1} =
-\frac{1}{n}\sum\strut^\prime v_k^2 + (2+\alpha)\sum_{k=1}^n v_k^2\,,
$$
where $\sum\strut^\prime$ denotes the summation over those $k$, for which
the $k$-th row of the matrix  $\tilde{\Cal D}_n$ contains
$\frac{1}{n}$ in the diagonal (recall that there are  $\sqrt n/3$ such~$k$).
Rewrite the last equality in the form
$$
\sum\limits_{k=1}^n (v_k-v_{k+1})^2 + \alpha\sum\limits_{k=1}^n v_k^2 =
\frac{1}{n}\sum\strut^\prime v_k^2\,. \tag 7
$$
An upper bound for the r.h.s.\ is evident:
$$
\frac{1}{n}\sum\strut^\prime v_k^2 \leq \frac{1}{n}\cdot\frac{\sqrt n}{3}=
\frac{1}{3\sqrt n}\,. \tag 8
$$
Using, if necessary, the cyclic permutation of the coordinates, we can think
that $v_1=\max\limits_{k}v_k$.
Then for any $s$ we have the following lower bound for the l.h.s.:
$$
\sum\limits_{k=1}^n (v_k-v_{k+1})^2 + \alpha\sum\limits_{k=1}^n v_k^2
\geq \sum\limits_{k=1}^s (v_k-v_{k+1})^2 \geq
\frac{(1-v_s)^2}{s}\,. \tag 9
$$
The second inequality here is the inequality about the arithmetic mean
and the quadratic mean (or its simple consequence, if some of the differences
$v_k-v_{k+1}$ are negative). Combining (7), (8) and (9), we obtain
inequality
$$
(1-v_s)^2\leq\frac{s}{3\sqrt n}\,,
$$
hence, taking into account that $0\leq v_s\leq 1$,
$$
v_s\geq 1-\sqrt{\frac{s}{3\sqrt n}}\,.
$$
As in the proof of the previous lemma we define
$v_k$ for all  $k\in {\Bbb Z}$ by the periodicity.
Then for  $|s|\leq \frac{5}{4}\sqrt n$
$$
v_s\geq  1-\sqrt{\frac{|s|}{3\sqrt n} }\geq 1-\sqrt{5/12}\,.
$$
Now we can give one more lower bound for the equality (7):
$$
\frac{1}{3\sqrt n}\geq
\frac{1}{n}\sum\strut^\prime v_k^2 =
\sum\limits_{k=1}^n (v_k-v_{k+1})^2 + \alpha\sum\limits_{k=1}^n v_k^2 \geq
\alpha \sum\limits_{k=-\frac{5}{4}\sqrt n}^{\frac{5}{4}\sqrt n} v_k^2 \geq
\Bigl(1-\sqrt{5/12}\,\Bigr)\alpha\cdot \frac{5}2\sqrt n\,.
$$
Thus,
$$
\alpha \leq \frac{2}{15(1-\sqrt {5/12})}\cdot\frac{1}{n}\approx
0.376\cdot\frac{1}{n}< \frac{2}{5n}\,.
$$
\hbox{}\hfill $\square$ Lemma 2.3 is proven.
\enddemo
 
\bigskip
\subhead
3. Combinatorial realization of the characteristic polynomial
\endsubhead

In this section we give a combinatorial realization of the characteristic
polynomial $P_{n, q}$ of the matrices $T_q(\tilde{\Delta}_n)$.

Consider graph which is a cycle with  $n$ vertices.
Enumerate its vertices along the cycle, assign the weight
$x-2\cos\frac{2\pi q k}{n}$ to the $k$-th vertex and weight $-1$ to each edge.
We denote this labeled graph $\Gamma_{n, q}$. Assign to an arbitrary
set $\xi$ of several non-intersecting edges of  $\Gamma_{n, q}$
the weight $w_\xi$, defined as a product of weights of all edges, belonging to
this set, and weights of all vertices not belonging to these edges.
For example, a graph $\Gamma_{5, 1}$ is drawn on the figure~2. The weight
of the set, consisting of two bold edges equals to $x-2\cos\frac{8\pi}{5}$,
and the weight of the empty set is $\prod_{k=1}^5 (x-2\cos\frac{2\pi k}{5})$.
Denote $K_{n, q}(x)=\sum_\xi w_\xi$, where the summation runs over all
possible sets of non-intersecting edges.

\proclaim{Lemma 3.1} $P_{n, q}(x) = K_{n, q}(x)+2(-1)^n$.
\endproclaim

\demo{Proof} Denote for shortness
$x{\Cal E}_n-T_q(\Delta_n)=(a_{ij})$. Then by the definition of determinant
$$
P_{n, q} = \sum_{\sigma\in S_n}
(-1)^\sigma a_{1\sigma(1)}a_{2\sigma(2)}\dots a_{n\sigma(n)}\,.\tag {10}
$$
Recall that the matrix $(a_{ij})$ is tridiagonal.
Consider an arbitrary nonzero monomial from the r.h.s.\ of this equality.
Assume that it contains diagonal factor $a_{kk}$ and off-diagonal factor
$a_{k+1\,k+2}$. Then it necessary contains the off-diagonal factor
$a_{k+2\,k+1}$. It is easy to see that all off-diagonal entries of matrix
$(a_{ij})$ enter every nonzero monomial in pairs of the form
$a_{k\,k+1}$, $a_{k+1\,k}$ (excluding two monomials:
$a_{12}a_{23}\dots a_{n1}$ and $a_{21}a_{32}\dots a_{1n}$).
This observation allows us to construct a bijection between the sets of
non-intersecting edges of our cycle $\Gamma_{n, q}$ and the nonzero monomials
in the equality~(10). Namely, the set
$(k_1, k_1+1)$, $(k_2, k_2+1)$, \dots corresponds to the monomial
$a_{11}\dots a_{k_1-1\,k_1-1}a_{k_1\,k_1+1}a_{k_1+1\, k_1}
a_{k_1+2\,k_1+2}\dots a_{k_2,k_2+1}a_{k_2+1\,k_2}\dots$.
The sign of this monomial is determined by the parity of  number of edges
and equals to the corresponding weight. The summand $2(-1)^n$
in equality (10) corresponds to the two exclusive monomials.
\hfill $\square$~\hbox{Lemma 3.1 is proven.}
\enddemo

\bigskip
\subhead
4. An estimation of spectral measure in the case of the infinite Heisenberg
group
\endsubhead

Let $\delta_e$ be the characteristic function of the group unity,
viewed as an element of the space $L^2$ over this group.
We will use this notation for different groups since it does not lead to a
confusion. In this section we estimate the value $(E_t\delta_e,\delta_e)$,
where $E_t$ is a spectral projector of the operator $\Delta$,
corresponding to the set  $[-1, -1+t)\cup (1-t, 1]$.

Remark that for the Laplace operator in ${\Bbb Z}^n$
built from  the standard generators it is easy to demonstrate with the help
of Fourier transform that $(E_t\delta_e,\delta_e)\sim \hbox{const }t^{n/2}$.

The next simple lemma is a key element for a transition from the infinite
group to a finite one.

\proclaim{Lemma 4.1} Let  $n$ be arbitrary positive integer, $N=n^2+1$.
Then $(\Delta^k\delta_e,\delta_e)=(\Delta_{N}^k\delta_e,\delta_e)$
for all $k\leq n$; the l.h.s.\ here is a scalar product in the space
$L^2(H)$, and the r.h.s.~--- in $L^2(H_N)$.
\endproclaim

\demo{Proof} The operator $\Delta^k$ regarded as an element of the group
algebra of the group $H$ is a formal linear combination of matrices of the
form
$\Bigl(\smallmatrix
1&a&c\cr
0&1&b\cr
0&0&1
\endsmallmatrix\Bigr)$,
where each matrix is a product of  $k$ generators or their inverse.
The scalar product, we are interested in, equals to the coefficient of the
identity matrix in such a decomposition.
It is clear that the entries  $a$ and $b$ in every such a matrix do not
exceed~$k$, and the entry $c$ does not exceed certain quadratic function of~$k$.
It is easy to prove by induction that $|c|\leq k^2/2$.
The difference between calculations for $L^2(H)$ and $L^2(H_N)$ is that
the entries of matrices in $H_N$ are not numbers but residues
in the range from  $-N/2$ to $N/2$.
Therefore for $k\leq n$ the results for $H$ and $H_N$ coincide.
\hbox{}\hfill $\square$ Lemma 4.1 is proven.
\enddemo

%\vfill\eject

\proclaim{Lemma 4.2} For any $\alpha>0$ there exist a constant
$c_0>0$, arbitrary large positive integer $n$, and a polynomial $P_n(x)$
of degree~$n$, such that have the following properties:

\item{1)} $|P_n(x)|\leq 1$ for $|x|\leq 1$, $P_n(1)=1, P_n(-1)=1$;
\item{2)} $|P_n(x)|\leq \frac{1}{n^6}$ for
          $|x|\leq 1-\frac{1}{n^{2-\alpha}}$;
\item{3)} $P_n(x)\geq с_0$ for $1-\frac{1}{n^{2}} \leq |x|\leq 1$.
\endproclaim

\demo{Proof}
%Without loss of generality we will lose some generality.
We chose $P_n$ as a properly normalized Chebyshev polynomial.
Namely,
$P_n(x)=С_n T_n\bigl(\frac{n^{2-\alpha}}{n^{2-\alpha}-1}\,x\bigr)$,
where $C_n$ is a normalizing constant.
For $|x|\geq 1$ the Chebyshev polynomials are given by the formulae
$$
T_n(x)=\frac{1}{2}\Bigl(\bigl(x+\sqrt{x^2-1}\,\bigr)^n+
        \bigl(x-\sqrt{x^2-1}\,\bigr)^n\Bigr)\,.
$$
Hence for $x=\pm 1$ and $n\to+\infty$ we derive:
$$
P_n(x)\sim\frac{C_n}2 \biggl(1+\frac{\sqrt2}{n^{1-\alpha/2}}\biggr)^n
\sim C_n e^{\sqrt2 n^{\alpha/2}}\,. \tag {11}
$$
Thus, we can take $C_n\sim 2e^{-\sqrt2 n^{\alpha/2}}$.
This guarantees not a polynomial but an exponential decreasing of $P_n(x)$
for  $|x|<1-\frac{1}{n^{2-\alpha}}$ (more precisely  for
$|x|<1-\frac{2}{n^{2-\alpha}}$) and the separation from zero
for $1-\frac{1}{n^{2}} \leq |x|\leq 1$.
\hbox{}\nobreak\hfill\hbox{$\square$ Lemma 4.2 is proven}.
\enddemo

\demo{Remark}  From minimax properties of the Chebyshev polynomials it
follows that the exponent $2-\alpha$ in the condition~2 of lemma 4.2
can not be replaced with~2.
Indeed, in this case an estimation similar to (11) shows that $C_n$
does not tend to 0, and therefore  $P_n$, as well as any other polynomial
satisfying the condition~1, is not small inside the interval~$[-1,1]$,
even if the additional condition 3 is not imposed.
\enddemo

\proclaim{Theorem 4.3} For any $\alpha>0$ and sufficiently small~$t$ the
estimation
$c_1 t^{2+\alpha}\leq (E_t\delta_e,\delta_e)$ holds.
\endproclaim

\demo{Proof}
Denote $t=1/n^2$, $N=n^2+1$, let $\chi_t$ be the characteristic function
for the set $A_t=[-1,-1+t]\cup[1-t,1]$. The key idea of the proof
is to substitute the function $\chi_t$ by a polynomial and then
``move'' all computations from the infinite Heisenberg group to the finite one.

The polynomials $P_n$ from lemma 4.2 obey the inequalities
$$
с_0\chi_t(x)\leq P_n(x)\leq
\chi_{t^{1-\alpha}}(x)+\frac{1}{n^{6}}\,.
$$
Consequently
$$
с_0\mu(A_t) \leq \|P_n(\Delta)\|^2\leq \mu(A_{t^{1-\alpha}})+\frac{1}{n^{6}}\,.
$$
Moreover, the same inequality holds for the operator $\Delta_N$.
Since $\deg P_n = n$, then in the spirit  of lemma~4.1 we have
$$
\|P_n(\Delta)\|^2=\|P_n(\Delta_N)\|^2\,.
$$
Thus,
$$
\mu(A_{t^{1-\alpha}})+\frac{1}{n^{6}} \geq
\|P_n(\Delta)\|^2=\|P_n(\Delta_N)\|^2 \geq
с_0 \mu_N(A_t)\,. \tag{12}
$$
The r.h.s.\ of this inequality is evaluated in the regular representation
of the finite Heisenberg group~$H_N$.
Notice, that thanks to lemma 2.3 in order to compute the r.h.s.\ it is
sufficient to employ the information only about the one-dimensional
representations of the group~$H_N$ (recall that
$t=\frac{1}{n^2}\asymp\frac{1}{N}$)!
Thus,
$$
\mu_N(A_t)\geq
\frac{1}{N^3}\!\!\!\!\!\!\sum_{s,t:\frac{s^2+t^2}{N^2}\leq
\frac{1}{n^2}}\!\!\!\! 1 \ \ \geq \frac{1}{n^4}\,.
$$
In the last inequality we estimate the number of terms as $N^2/n^2\asymp n^2$.
Together with (12) it leads to the required inequality.
Formally, in order to be able to apply lemma 2.3 we must take $N$ to be a
prime number. Therefore we had to choose $N$ more accurately, e.g.,
$N$ could be chosen as a prime number in the interval $[n^2+1, 2n^2+1]$.
\hbox{}\nobreak\hfill $\square$ Theorem 4.3 is proven.

\enddemo

Unfortunately one cannot obtain the estimation
$(E_t\delta_e,\delta_e) \leq c_2 t^{2-\alpha}$ in the same manner.
In order to avoid multi-dimensional representations, we need to take
too small~$t$. But according to the remark after lemma 4.2, the polynomial
$P_n$ in this case does not approximate the function $\chi_t$.

In the paper [K] the author proves the following inequality
for the Laplace operator on an amenable group:
$$
(E_t\delta_e,\delta_e)\geq
\frac{1-\frac{2\varepsilon}{t^2}}{|A_\varepsilon|}\,,
$$
where  $\varepsilon>0$ is any number, $|A_\varepsilon|$ is the size
of Folner's $\varepsilon$-set. The best value of the r.h.s.\
for fixed $t$ is achieved for $\varepsilon=O(t^2)$. In this case the
numerator is a certain constant and the denominator is presumably
$O(1/\varepsilon^4)$. This yields the estimation
$(E_t\delta_e,\delta_e)\geq Ct^8$, which is worse than our result.

\bigskip
\baselineskip=10pt
\rm \noindent Department of mathematics and mechanics,

\rm \noindent St.-Peterburg State University,

\rm \noindent St.-Petersburg, Russia.

\smallskip
\rm \noindent kostik\@kk1437.spb.edu

\end